\documentclass[letterpaper, 10 pt, journal,twoside]{IEEEtran}  

\IEEEoverridecommandlockouts                              

\usepackage{pstool}

\usepackage{amsmath}
\usepackage{amssymb}
\usepackage{algorithmic}
\usepackage{algorithm}
\usepackage{datetime}
\usepackage{acronym}
\newdateformat{mydate}{\THEDAY~\monthname[\THEMONTH]~\THEYEAR}
\mydate
\usepackage{multirow}
\usepackage{nth} 

\usepackage{caption}
\usepackage{subcaption}
\usepackage{epstopdf}
\usepackage{graphicx}

\usepackage{pgfplots}
\usepackage{tikz}
\usetikzlibrary{calc,arrows}
\usepackage{tikzscale}

\usepackage{ctable} 
\usepackage{todonotes} 
\usepackage[noadjust]{cite}	

\newcommand{\re}{\textrm{cmd}}
\newcommand{\8}{\textrm{ref}}
\newcommand{\no}{0}

\renewcommand{\baselinestretch}{1}

\title{Predictive Control of Autonomous Kites in Tow Test Experiments${}^\star$}

\author{Tony A. Wood$^{1}$, Henrik Hesse$^{2}$, and Roy S. Smith$^{1}$
\thanks{${}^\star$This manuscript is a preprint of a paper accepted for the IEEE Control Systems Letters and is subject to IEEE Control Systems Society copyright. 
Upon publication, the copy of record will be available at http://ieeexplore.ieee.org.}
\thanks{$^{1}$T. A. Wood and R. S. Smith are with the Automatic Control Laboratory (IfA), ETH Zurich, Physikstrasse 3, 8092 Zurich, Switzerland. E-mail: {\{woodt, rsmith\}@control.ee.ethz.ch}. Tel: +41 44 632 3629.}
\thanks{$^{2}$H. Hesse is with the Aerospace Sciences Division, University of Glasgow (Singapore), 510 Dover Road, Singapore 139660 (formerly at IfA, ETH Zurich). Email: henrik.hesse@glasgow.ac.uk}
\thanks{This research was supported by the Swiss National Science Foundation (Synergia) No. 141836.
}}

\begin{document}
\maketitle
\thispagestyle{empty}
\pagestyle{empty}

\acrodef{mpc}[MPC]{Model Predictive Control}
\acrodef{gs}[GS]{ground station}
\acrodef{fhnw}[FHNW]{Fachhochschule Nordwestschweiz}

\begin{abstract}
In this paper we present a model-based control approach for autonomous flight of kites for wind power generation. Predictive models are considered to compensate for delay in the kite dynamics. We apply Model Predictive Control (MPC), with the objective of guiding the kite to follow a figure-of-eight trajectory, in the outer loop of a two level control cascade. The tracking capabilities of the inner-loop controller depend on the operating conditions and are assessed via a frequency domain robustness analysis. We take the limitations of the inner tracking controller into account by encoding them as optimisation constraints in the outer MPC. The method is validated on a kite system in tow test experiments.
\end{abstract}


\acresetall 


\section{Introduction}
Airborne Wind Energy (AWE) generators have been proposed as a mobile, cost-effective, and more sustainable alternative to conventional wind turbines. In this work we focus on kite power systems which generate power by flying a multi-line tethered wing or kite in crosswind motion following a figure-of-eight pattern. The tethers are connected to winches on the ground which generate power by unreeling the lines in the so-called traction phase. The traction phase is alternated with a retraction phase to reel in the tethers using only a fraction of the energy generated during traction, leading to a net positive cycle power. 

Since the original inception of the kite power concept~\cite{Loyd80}, several groups have developed prototype systems 
and the reference book~\cite{AWE13} provides a detailed overview of the field of~(AWE). In this work we focus on the control of kites during the traction phase of a ground-based system, as developed in~\cite{A2WE}, with ground-based measurements of line angles and length, steering actuation, and power generation at the so-called \ac{gs}. The experimental implementation of autonomous kite power systems requires stabilising control approaches which can handle the unstable, nonlinear dynamics. 
Model-free guidance approaches based on a switching point strategy and using the kite heading angle as feedback variable for tracking~\cite{FagEtal:2014:IFA_4502, Erhard201513} provide a successful starting point for further control development. 
The performance of such model-free approaches, however, is affected by varying operating conditions and time delay 
which make tuning difficult and effectively limit the overall system power output, especially for ground-based systems. 

To improve the control performance in experimental implementations of kite power prototypes, model-based control approaches have been proposed in \cite{Erhard201513,woodACC15,woodCDC15,RontsisEtAlCDC15}. The underlying kinematic models link the kite heading angle to the overall kite motion. In \cite{woodACC15} the estimated 
kite heading angle was further related to the steering input with a model that includes an input delay. 
State estimation for ground-based (AWE) system with output delay is addressed in \cite{polzinIFACWC17}.
In~\cite{CostelloECC2015,woodACC15,woodCDC15,RontsisEtAlCDC15,woodIFACWC17} guidance strategies have been developed for experimental kite power systems to allow path following. 
In particular, \cite{woodCDC15} used a kinematic model including input delay for a figure-of-eight path planning and tracking strategy. To ensure robust performance of the control approach we further considered limitations on the tracking bandwidth imposed by the input delay. 
The robustness of the cascaded control architecture of~\cite{woodACC15,woodCDC15} was further improved following a \ac{mpc} approach in~\cite{woodIFACWC17} by formulating the path following problem as an optimisation problem with constraints on the heading angle, its rate of change, and the kite position. 

\ac{mpc} is naturally suited to address the complexities in the control of kites, e.g. constraint satisfaction and minimisation of the deviation from a reference path, and has been extensively explored in simulation in~\cite{diehl2005real,canale2010high,ilzhofer2007nonlinear,Zanon2013}. 
The application of \ac{mpc} approaches to kite power systems, however, tends to be sensitive to unmodelled dynamics and hindered by limitations in processing power for real-time operation. In \cite{woodIFACWC17} we therefore use the kinematic model introduced in~\cite{woodCDC15} which allows for online adaptation of model parameters to reduce the model mismatch. The approach was demonstrated in simulation to achieve path following while satisfying constraints imposed by the limitations of a lower-level tracking controller that are subject to model uncertainty and input delay. The ability to account for constraints is the main benefit over computationally more efficient guidance methods that are purely based on geometry such as in~\cite{jehle2014applied}.

In this paper we extend the kinematic model to account for variations of the kite velocity within a figure-of-eight cycle which significantly improves the tracking capabilities of the path following controller developed in \cite{woodIFACWC17}. Moreover, as main contribution, we demonstrate the performance of the \ac{mpc} approach for varying line length during tow test experiments with a prototype kite power system. 

In Section~\ref{sec:controlScheme} we describe the cascaded control architecture with delay compensation and constrained outer-loop guidance that accounts for limitations of the inner-loop controller. In Section~\ref{sec:experimentalImplementation} we describe the implementation of the control scheme for tow test experiments. We present experimental results in Section~\ref{sec:results} before concluding in Section~\ref{sec:conclusion}.
    
\section{Control Scheme}\label{sec:controlScheme}
\begin{figure}[tb]
  \centering
  \includegraphics[width=0.72\linewidth]{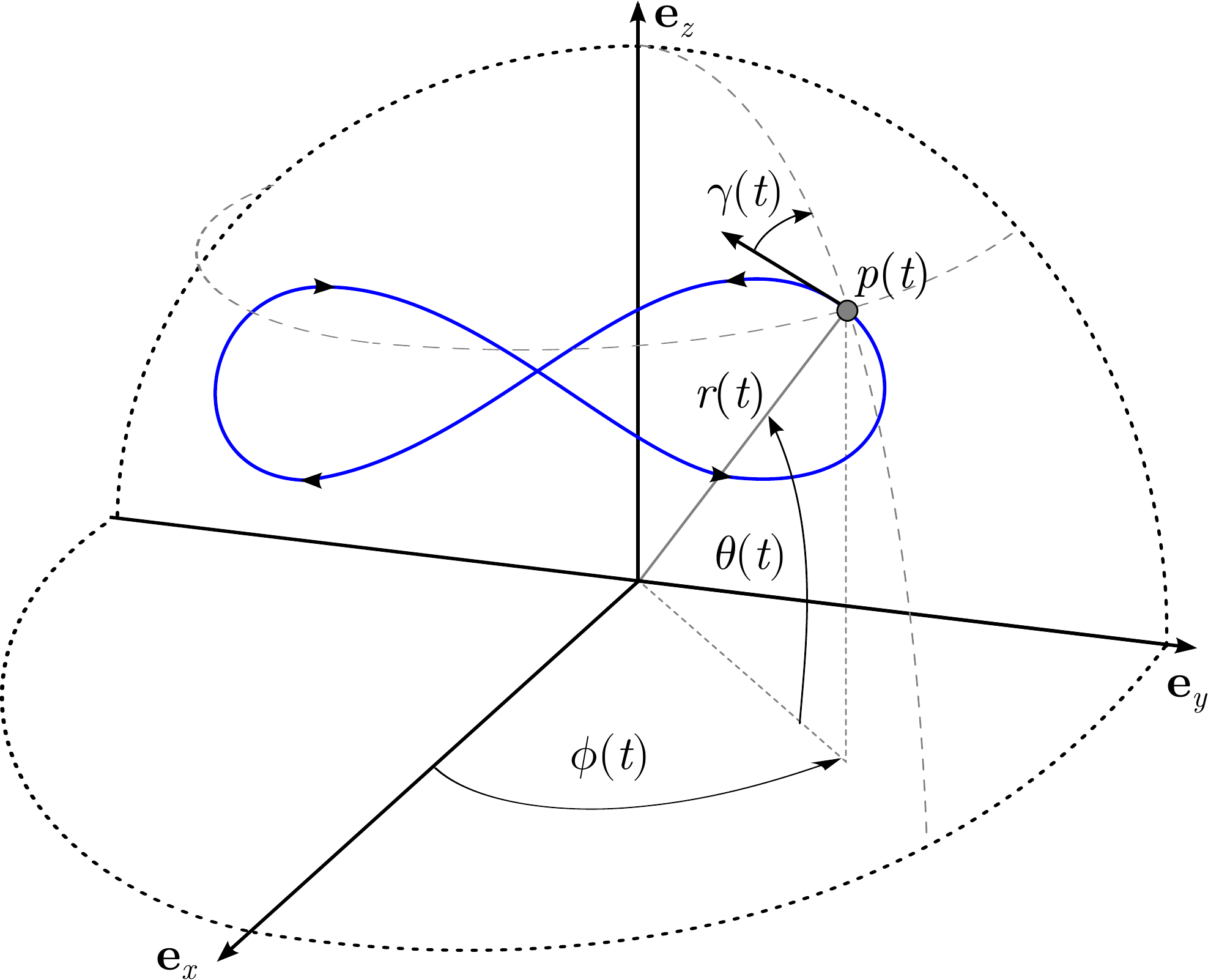}
  \caption{Coordinate system, with $x$-axis aligned with the wind direction, and kite trajectory showing kite position, $p(t)$, and heading angle, $\gamma(t)$.\vspace{0pt}}
  \label{fig:coordinates}
\end{figure}

We consider a two line kite power system in the traction phase where the kite is flown in crosswind conditions. The motion of the kite perpendicular to the tethers is actuated by the difference in lengths of the two lines, $\delta(t)$. The position of the kite in Cartesian coordinates, $p(t)$, can be expressed as a function of the elevation angle, $\theta(t)$, the azimuth angle, $\phi(t)$, and the line length, $r(t)$, which are illustrated in Figure~\ref{fig:coordinates} and are all measured from the \ac{gs}, 
 \begin{align*}
    p(t) =
    \begin{bmatrix}
      p_x(t)\\
      p_y(t)\\
      p_z(t)
    \end{bmatrix} =
    \begin{bmatrix}
      r(t)\cos(\theta(t))\cos(\phi(t))\\
      r(t)\cos(\theta(t))\sin(\phi(t))\\
      r(t)\sin(\theta(t))
    \end{bmatrix}\,.
  \end{align*}

\begin{figure}[tb]
  \centering
  \includegraphics[width=\linewidth]{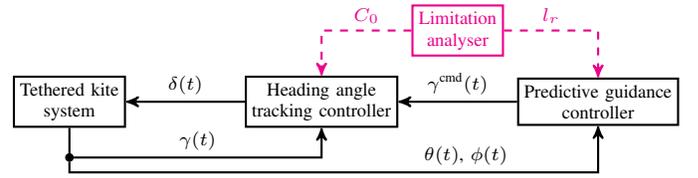}
  \caption{Cascaded control architecture with inner-loop controller tracking a commanded heading angle, $\gamma^\re(t)$, given by the outer-loop guidance controller which incorporates the limitations of the inner loop as constraints, $|\dot{\gamma}^\re(t)|<l_r$.\vspace{0pt}}
  \label{fig:cascade}
\end{figure}

A cascaded architecture, illustrated in Figure~\ref{fig:cascade}, is used to control 
the kite along a figure-of-eight path. Controlling the heading angle 
defined as
\begin{align}\label{eq:gamma}
  \gamma(t) := \arctan\left(\frac{\cos(\theta(t))\dot{\phi}(t)}{\dot{\theta}(t)}\right)\,,
\end{align} has been shown to be an effective approach for autonomous crosswind flight control of kites \cite{FagEtal:2014:IFA_4502,Erhard201513,woodACC15}. The outer loop of the cascade is controlled by the guidance controller which produces a commanded heading angle trajectory, $\gamma^\re(t)$, that is tracked in the inner loop by the tracking controller. 

The steering behaviour of the kite is affected by line dynamics due to the indirect tether actuation. This effect can be modelled as an input delay and taken into account in a model-based delay compensation scheme~\cite{woodACC15}. The presence of delay and model uncertainty imposes fundamental limitations on the tracking performance. As in~\cite{woodIFACWC17} we assess the limitations of the tracking controller, based on estimates of the current operating conditions, 
with a frequency domain robustness analysis. We parametrise the limitations by an upper bound on the rate of change of the commanded signal. 
The rate limit, $l_r$, is communicated to the optimisation-based 
guidance controller which takes the current inner-loop 
tracking capabilities into account as constraints. 
 
\subsection{Control Models}
The evolution of the kite position expressed in terms of line angles, $\xi(t):=(\theta(t),\phi(t))$, can be modelled by the following kinematic unicycle model, \cite{woodCDC15},
\begin{subequations}\label{eq:unicycle}
  \begin{align}
  \dot{\theta}(t) &= \frac{v_{\theta\phi}(t)}{r}\cos\left(\gamma(t)\right)\,,\\
  \dot{\phi}(t) &= \frac{v_{\theta\phi}(t)}{r\cos(\theta(t))}\sin\left(\gamma(t)\right)\,,
\end{align}
\end{subequations} where the velocity perpendicular to the tethers is modelled as a static function of the position, $\xi(t)$, and orientation, $\gamma(t)$,
\begin{align}\label{eq:velocity}
  v_{\theta\phi}(t) = r\alpha_L \cos\left(\theta(t)\right)\cos\left(\phi(t)\right) - r\alpha_G\cos\left(\gamma(t)\right)\,, 
\end{align} with parameters $\alpha_L,\alpha_G>0$ representing velocity components arising from lift and gravitational forces.
Modelling the velocity to vary within a figure-of-eight cycle captures the true behaviour more closely than the assumption of constant velocity made in the control model of~\cite{woodCDC15, woodIFACWC17} as can be seen in the experimental results shown in Figure~\ref{fig:velocity1C}. 

Based on experimental observations in \cite{woodACC15}, we relate the heading angle, $\gamma(t)$, to the steering input, $\delta(t)$, by integrator dynamics with a time delay $t_d$,
\begin{align}
  \label{eq:steeringModel}
    \dot{\gamma}(t) &= K\delta(t-t_d)\,.
\end{align} 
The 
parameters, $r$, $\alpha_L$, $\alpha_G$, $K$, and $t_d$ are 
considered constant within the control horizon and are re-identified online based on updated measurements. 

\subsection{Tracking Controller}
To control the heading angle we account for the delay in \eqref{eq:steeringModel} in a model-based approach as described in~\cite{woodACC15}. We predict the orientation after the delay time, $\gamma^{t_d}(t)$, and apply a proportional gain controller to the difference of the commanded and the predicted output,
\begin{align}\label{eq:smithPredictorRealisation}  
  \delta(t)=C_0(\gamma^\re(t)-\gamma^{t_d}(t))\,.
\end{align}
Controlling delayed systems based on model predictions is referred to as predictor feedback \cite{krstic2009delay}. To assess the performance of the tracking controller we consider the regulation of the time-shifted tracking error, $e_{t_d}(t) := \gamma^\re(t)-\gamma(t+t_d)$.

The control model in~\eqref{eq:steeringModel} is simple but there is a considerable degree of uncertainty in the steering gain parameter, $K$, and delay, $t_d$. 
In the frequency domain we model the uncertainty of the plant, $\gamma=G(s)\delta$, 
\begin{align*}
  G(s)&=G_\no(s)e^{-st_d}=\frac{K}{s}e^{-st_d}\,,
\end{align*} 
 with a 
multiplicative perturbation with the perturbation weight for delayed first-order systems introduced in \cite{laughlin1987smith},
\begin{align}
  \mathcal{G} &\!:=\! \left\{\left(1\!+\!W_m(s)\Delta(s)\right)G(s)\Big| \|\Delta(s)\|_\infty\leq 1\right\},\label{eq:uncertaintySet}\\
  W_m(j\omega) &\!\:=\!
  \begin{cases}
    \left|\frac{K+\delta K}{K}e^{-j\delta t_d \omega}-1\right|&\textrm{if }\omega< \frac{\pi}{\delta t_d},\\
    \left|\frac{K+\delta K}{K}\right|+1&\textrm{if }\omega\geq\frac{\pi}{\delta t_d},
  \end{cases}\nonumber
\end{align} where $\delta K$, $\delta t_d$ are the bounds on the deviations of the parameters, and $\Delta(s)$ is an unknown but bounded perturbation. 

The 
controller in~\eqref{eq:smithPredictorRealisation}, parametrised by the control gain, $C_0$, 
can be written 
as, $\delta = C(s;C_0) (\gamma^\re-\gamma)$, with
\begin{align*}
  C(s;C_0) &= \frac{C_0}{1+C_0 G_0(s)(1-e^{-st_d})}\,.
\end{align*} For 
perturbed plants, $G^p(s)\in\mathcal{G}$, the 
relation 
between the input and the shifted 
error is, $e_{t_d}=S^p_{t_d}(s;C_0)\delta$, with
\begin{align*}
  S^p_{t_d}(s;C_0)&=\frac{1+L(s;C_0)(1-e^{st_d})\left(1+W_m(s)\Delta(s)\right)}{1+L(s;C_0)\left(1+W_m(s)\Delta(s)\right)}\,,
\end{align*}where $L(s;C_0)\!=\!C(s;C_0)G(s)$ is the loop transfer function.

The tracking capabilities depend on the control gain, $C_0$, and on the properties of the commanded signal, $\gamma^\re(t)$. We would like to set the control gain and constrain the commanded signal such that the shifted error remains small, 
$|e_{t_d}(t)|<l_e$, for all perturbed plants in the set given in~\eqref{eq:uncertaintySet}. In particular, we limit the magnitudes of the commanded signal, $|\gamma^\re(t)|<l_m$, and its rate of change, $|\dot{\gamma}^\re(t)|<l_r$. As discussed in~\cite{woodIFACWC17}, these limits on the commanded signal and the shifted tracking error 
are approximately translated to the frequency domain with the robust performance condition
\begin{align}\label{eq:RP}
  \sup_{\|\Delta\|_\infty\leq 1} \|W_p(j\omega;l_r)S_{t_d}^p(j\omega;C_0)\|_\infty&< 1\,,
\end{align}  for specified parameters $l_r$ and $C_0$, where we design the performance weight to be
\begin{align*}
    W_p(j\omega;l_r) &=
  \begin{cases}
    \frac{l_m}{l_e} &\textrm{if }\omega< \frac{l_r}{l_m}\,,\\
    \frac{l_r}{l_e\omega}&\textrm{if }\omega\geq \frac{l_r}{l_m}\,.
  \end{cases}
\end{align*} 

Bounds for the left-hand side of~\eqref{eq:RP} can be determined by a structured singular value analysis. 
Given values of the model parameters, $K$, $t_d$, the level of uncertainty considered, $\delta K,\delta t_d$, and the bounds $l_m$, $l_e$, the condition in \eqref{eq:RP} depends on the control gain, $C_0$, and the rate limit, $l_r$. We follow the scheme suggested in~\cite{woodIFACWC17}, where $C_0$ is chosen for robust stability and $l_r$ is selected such that robust tracking performance is guaranteed, i.e. the condition in~\eqref{eq:RP} is satisfied. 

\subsection{Predictive Guidance Control} 
For the design of the guidance controller we assume that the tracking controller is able to follow the commanded signal well but delayed by the steering delay, i.e., $\gamma(t)\approx\gamma^\re(t-t_d)$. We apply the model given in~\eqref{eq:unicycle} and \eqref{eq:velocity} to compensate for the delayed tracking. Given 
a prediction of the kite position $t_d$ ahead of time, $\xi^{t_d}(t):=(\theta^{t_d}(t),\phi^{t_d}(t))$, we find a value for the commanded signal at the current time, $\gamma^\re(t)$, that controls the kite to follow a reference figure-of-eight path, $\xi^\8:=(\theta^\8,\phi^\8)$. The reference path is generated and updated online such that the corresponding reference orientation, $\gamma^\8$, is of sinusoidal form and its rate of change satisfies the limitation given by the tracking controller as suggested in~\cite{woodCDC15}.  
This implies that for lower limits on the rate of change, larger figure-of-eight paths are required.

We apply a \ac{mpc} approach with the objective of minimising the deviation of the delay compensating prediction, $\xi^{t_d}(t)$,  from the reference path, $\xi^\8$, over a finite time horizon, $T_H$, while constraining  the kite to remain in a predefined safety window, $\underline{\xi}\leq\xi^{t_d}(t)\leq\overline{\xi}$, and satisfying the limit on the rate of change of the commanded heading angle imposed by the tracking controller, $|\dot{\gamma}^\re(t)|<l_r$. We also constrain the magnitude of the commanded orientation, $|\gamma^\re(t)|<l_m$, to avoid commanding the kite to fly straight down towards the ground. 
More details on the optimisation problem that is formulated and solved to capture the guidance objective are presented in Section~\ref{sec:controllerImplementation} where the kinematic model given in~\eqref{eq:unicycle} and~\eqref{eq:velocity} is used in the formulation of state constraints and the limitation on the turn rate, derived from steering model~\eqref{eq:steeringModel}, is implemented as input rate constraint. 

\begin{figure}[tb]
  \centering
  \includegraphics[width=\linewidth]{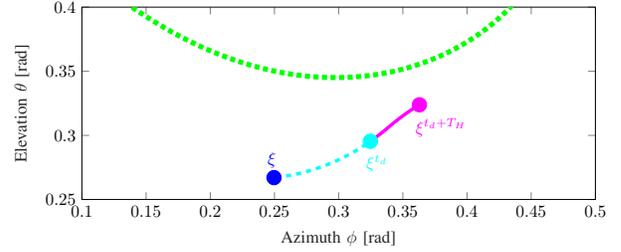}
  \caption{Prediction of flight trajectory: delay compensation ({\color{cyan}dashed}),   \ac{mpc} prediction ({\color{magenta}solid}); reference path ({\color{green}dotted}).\vspace{0pt}}
  \label{fig:guidancePath}
\end{figure}

Figure~\ref{fig:guidancePath} illustrates the guidance concept with the predicted path of the kite. 
The prediction consists of two parts: 
the first part is determined by the past inputs in the delay compensation scheme; 
the further evolution is due to current and future actuation and can be optimised by the \ac{mpc}.

\section{Experimental Implementation}\label{sec:experimentalImplementation}
 \ac{mpc} for autonomous kites has been studied extensively in literature \cite{diehl2005real,canale2010high,ilzhofer2007nonlinear,Zanon2013} but demonstrations have been limited to simulation. In this work we 
design and conduct experiments to validate the applicability of predictive crosswind flight control. 
We further demonstrate a novel experimental method to test kite controllers in low wind conditions using a tow test configuration. 

\subsection{Tow Test Experiments}
\begin{figure}[tb]
  {    \centering
  \hfill
  \includegraphics[height=0.39\linewidth]{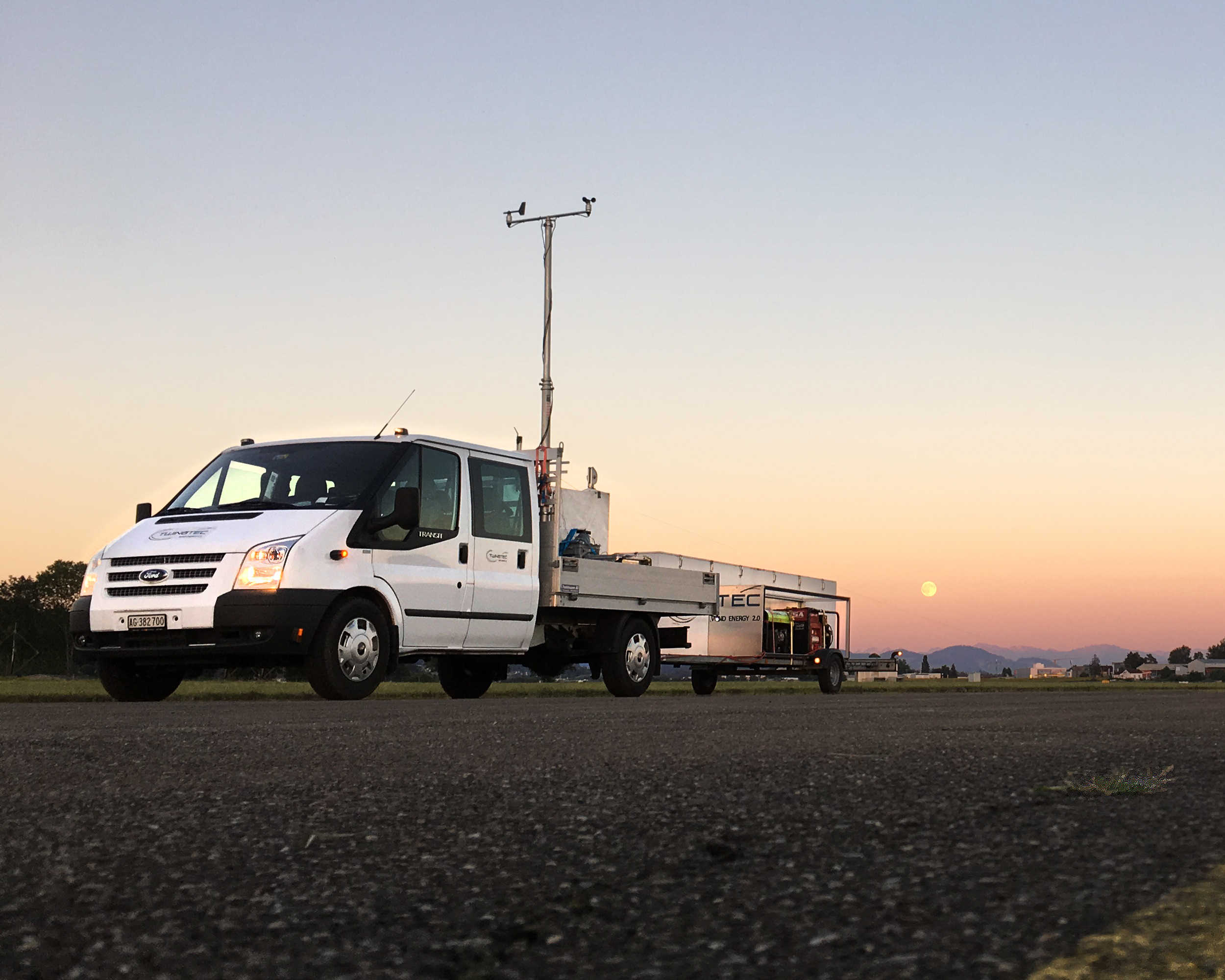}
  \includegraphics[height=0.39\linewidth]{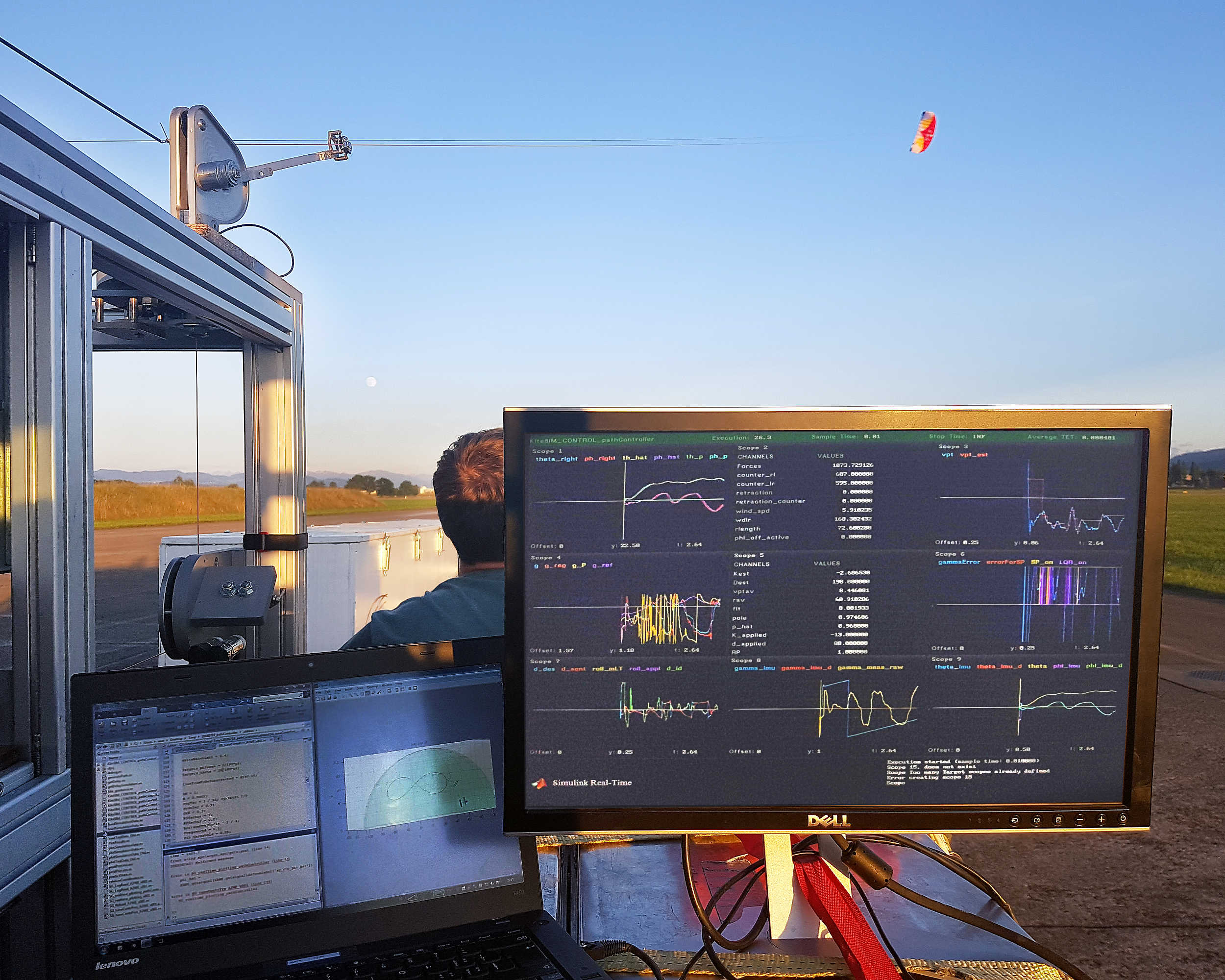}
  \hfill}
  \caption{Tow test configuration with the kite connected to \ac{gs} control system (right) mounted on a moving vehicle (left).\vspace{0pt}}
\label{fig:towTest}
\end{figure}

The experimental flight control tests are performed on a prototype (AWE) system \cite{A2WE} which has been developed at \ac{fhnw} and modified to allow for tow test experiments as depicted in Figure~\ref{fig:towTest}. By mounting the \ac{gs} on the back of a truck we can create relative wind when driving the vehicle on an airfield runway. The origin of the coordinate frame shown in Figure~\ref{fig:coordinates} moves with the truck and the $x$-axis points in the opposite direction to the vehicle velocity. 

For the results presented here, an HQ Apex III 5m${}^2$ ram air kite was connected to the \ac{gs} via two tethers of maximum length 150~m.  
The forces on the tethers are controlled by regulating the reel-out speed of the winches. The winch control is independent of the steering control. The entire control and estimation architecture 
is implemented in Matlab Simulink and runs on a Speedgoat Real-Time Target Machine with a fixed sampling time, $T=0.01$~s. With this mobile 
configuration we can test in still wind conditions and create reproducible relative wind scenarios to develop controllers for crosswind flight of kites with variable tether length.

\subsection{Controller Implementation}\label{sec:controllerImplementation}
As part of the predictive control approach in this work, we use real-time measurement data to identify and update the control model parameters online in an adaptive fashion. Given estimates of the model parameters, the tracking control gain, $C_0$, and the rate limit on the commanded heading angle, $l_r$, are set such that robust tracking performance is guaranteed. We use a sufficient analytic condition derived in~\cite{woodIFACWC17} to determine suitable values for $C_0$ and $l_r$. 
The condition involves an upper bound on the left-hand side of \eqref{eq:RP} which is explicitly parametrised by the tuning variables $C_0$ and $l_r$, allowing for simple evaluation online.



The kinematic model given in~\eqref{eq:unicycle} and~\eqref{eq:velocity} is discretised with the forward Euler method, 
\begin{align*}
      \theta_{k+1}\! &= \!\theta_k\! +\! T\alpha_L\cos(\theta_k)\cos(\phi_k)\cos(\gamma_k)\!-\!T\alpha_G\cos^2(\gamma_k)\,,\\
  \phi_{k+1}\! &= \!\phi_k\! +\! T\alpha_L\cos(\phi_k)\sin(\gamma_k)\!-\!\frac{T\alpha_G \sin(2\gamma_k)}{2\cos(\theta_k)}\,.
\end{align*}  With the short sampling period of $T=0.01$~s, no significant error is introduced by the discretisation. We denote the discretised model, which is used in the guidance control, as $\xi_{k+1} = f_d(\xi_k,\gamma_k)$, with $\xi_k := (\theta_k,\phi_k)$.

The figure-of-eight reference path is parametrised by a periodic sequence of positions, $\xi^\8=(\xi^\8_1,\xi^\8_2,\dots,\xi^\8_N)$, and orientations, $\gamma^\8=(\gamma^\8_1,\gamma^\8_2,\dots,\gamma^\8_N)$, that satisfy the model dynamics, $\xi^\8_{i+1}=f_d(\xi^\8_i,\gamma^\8_i)$, for $i=0,1,\dots,N-1$, and $\xi^\8_1 = f_d(\xi^\8_N,\gamma^\8_N)$.

To capture the guidance objective, we consider a system describing the deviation of the prediction of the position after the delay from the reference path with state $\chi_k:=\xi^{t_d}_k-\xi^\8_j$ and input $u_k:=\gamma^\re_k-\gamma^\8_j$, where $(\xi^\8_j,\gamma^\8_j)$ represent the reference point closest to the prediction of the kite state after the delay time. 
Minimising the deviation of the position prediction from the reference path corresponds to controlling the deviation system to zero. 

We linearise the kite kinematics around the reference path, starting at the reference point closest to the prediction of the kite position after the delay time, 
to obtain the linear time varying system $\chi_{i+1}\approx A_i\chi_i+B_iu_i$, with $A_i:=\frac{\partial f_d}{\partial \xi}(\xi^\8_{j+i},\gamma^\8_{j+i})$ and $B_i:=\frac{\partial f_d}{\partial \gamma}(\xi^\8_{j+i},\gamma^\8_{j+i})$, for $i=0,1,\dots,H-1$, where $H$ is the number of discrete time steps in the prediction horizon. To encode the input rate constraint into the optimisation problem we augment the deviation state such that its dynamics, $\hat{\chi}_{i+1}=\hat{A}_i\hat{\chi_i}+\hat{B}_i\Delta u_i$, 
\begin{align*}
\hat{\chi}_i&:=\begin{bmatrix} \chi_i\\u_{i-1}\end{bmatrix}\,,&
  \hat{A}_i&:=
        \begin{bmatrix}
          A_i & B_i\\
          0 & 1
        \end{bmatrix}\,,& \hat{B}_i&:=
  \begin{bmatrix}
    B_i\\
    1
  \end{bmatrix}\,,
\end{align*} are actuated by the change in input, $\Delta u_i:=u_i-u_{i-1}$.

The following optimisation solves the guidance task,
\begin{align}\label{eq:mpc}
   \min_{\Delta u,\epsilon}&&\!\!\sum_{i=0}^{H-1}\!&\big(\hat{\chi}_i^\top  \hat{Q}\hat{\chi}_i
\!+\!\epsilon_i^\top S\epsilon_i\big)\!+\!\hat{\chi}_H^\top \hat{Q}_H\hat{\chi}_H\!+\!\epsilon_H^\top S_H\epsilon_H,\\
   \operatorname{s.t.}&  &&\hat{\chi}_{i+1}=\hat{A}_i\hat{\chi}_i+\hat{B}_i\Delta u_i\,,\nonumber\\
   & &&\underline{\chi}_i \leq \chi_i + \epsilon_i\,, \quad \chi_i - \epsilon_i\leq\overline{\chi}_i \,,\quad \epsilon_i\geq 0\,,\nonumber\\
   &&&\underline{u}_i \leq u_i\leq \overline{u}_i\,,\quad\underline{\Delta u}_i \leq \Delta u_i\leq \overline{\Delta u}_i \nonumber\,,
\end{align} with actuation rate sequence $\Delta u := (\Delta u_0,\Delta u_1,\dots,\Delta_{H-1})$, slack variables sequence $\epsilon := (\epsilon_0,\epsilon_1,\dots,\epsilon_{H})$, state limits $\underline{\chi}_i:=\underline{\xi}-\xi^\8_{j+i}$, $\overline{\chi}_i:=\overline{\xi}-\xi^\8_{j+i}$, input limits $\underline{u}_i:=-l_m-\gamma^\8_{j+i}$, $\overline{u}_i:=l_m-\gamma^\8_{j+i}$, rate limits $\underline{\Delta u}_i:=-l_rT-\gamma^\8_{j+i}+\gamma^\8_{j+i-1}$, $\overline{\Delta u}_i:=l_rT-\gamma^\8_{j+i}+\gamma^\8_{j+i-1}$, and positive-definite weighting matrices $S$, $S_H$, $\hat{Q}:=\operatorname{diag}(Q,R)$, $\hat{Q}_H:=\operatorname{diag}(Q_H,R)$, 
where $(Q,Q_H)$ penalise the position deviation, $R$ penalises the orientation deviation, and $(S,S_H)$ penalise the slack variables. 
The inclusion of the slack variables 
with high penalisation $S\gg Q$ and $S_H\gg Q_H$ prevents the constraint optimisation problem from becoming infeasible in situations where there is no possibility of keeping the kite in the desired position window. This is relevant when initialising the controller from arbitrary positions. 
The optimisation problem \eqref{eq:mpc} is implemented and solved online with the optimisation software generation tool FORCES Pro \cite{FORCESPro}.

The value of the current commanded heading angle at time step $k$ is determined by the first element of the sequence ${\Delta u}^*$ that obtains the minimum of \eqref{eq:mpc}, $\gamma^\re_k={\Delta u}^*_0 + \gamma^\re_{k-1} - \gamma^\8_{j-1} + \gamma^\8_{j}$. The optimisation problem is re-solved in every time step with a receding horizon.

\section{Results}\label{sec:results}
The predictive control scheme has been experimentally tested in tow test experiments. In this section we present the results of an experiment conducted on the runway of St. Stephan airport in Switzerland on 9 December 2016. Data from two test flights (Flight~1 and Flight~2) are presented.

\begin{figure}[tb]
  \centering
  \includegraphics[width=\linewidth]{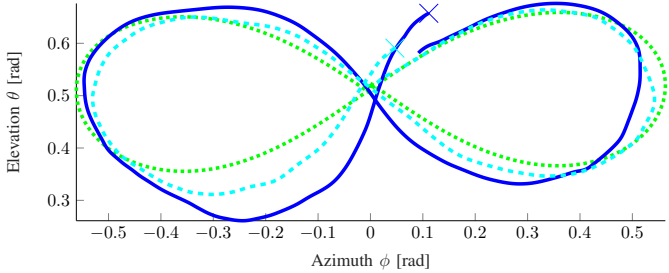}
  \caption{Trajectory tracking results (Flight~1): kite position trajectory, $\xi(t)$ ({\color{blue}solid}), following a reference figure-of-eight path, $\xi^\8$ ({\color{green}dotted}), with delay compensation prediction trajectory, $\xi^{t_d}(t)$ ({\color{cyan}dashed}). The trajectories start at the position marked with crosses at time $t=0$s.\vspace{0pt}}
 \label{fig:trajectory1C}
\end{figure}

\begin{figure}[tb]
  \centering
  \begin{subfigure}{\linewidth}
      \includegraphics[width=\linewidth]{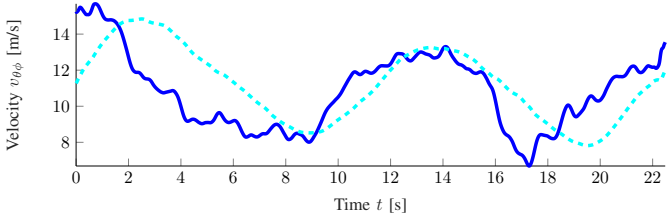}
  \caption{Velocity trajectory, $v_{\theta\phi}(t)$, estimated based on derivatives of position measurements, $r(t)\sqrt{\dot{\theta}^2(t)+\cos^2(\theta(t))\dot{\phi}^2(t)}$ ({\color{blue}solid}), and based on model~\eqref{eq:velocity} ({\color{cyan}dashed}).}
  \label{fig:velocity1C}
\end{subfigure}
    \begin{subfigure}{\linewidth}
      \vspace{15pt}
      \includegraphics[width=\linewidth]{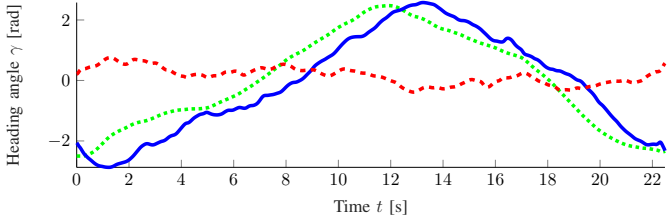}
  \caption{Trajectory of heading angle, $\gamma(t)$ ({\color{blue}solid}), tracking commanded signal, $\gamma^\re(t)$ ({\color{green}dotted}), with shifted tracking error, $e_{t_d}(t)$ ({\color{red}dashed}).}
  \label{fig:gamma1C}
  \end{subfigure}
  \caption{Trajectories of velocity and heading angle (Flight~1).\vspace{0pt}} 
  \label{fig:velocityOrientation1C}
\end{figure}

The number of stages in optimisation problem \eqref{eq:mpc} was selected to be $H=30$, resulting in a \ac{mpc} prediction horizon of $T_H=0.3$~s. Note that the overall prediction horizon of the guidance controller consists of the sum of the delay compensation and the \ac{mpc} horizon; with an average estimated delay of approximately $t_d=0.7$~s the overall prediction horizon stretches over approximately 1~s. The weights of the \ac{mpc} objective were $Q=\operatorname{diag}(1,2)$, $Q_H=5\cdot Q$, $R=5\cdot 10^{-3}$, $S=10^5\cdot Q$, $S_H = 10^5\cdot Q_H$. The safety window was defined by the position bounds $\underline{\theta}=0.17$~rad, $\overline{\theta}=1.40$~rad, $\underline{\phi}=-0.70$~rad, and $\overline{\phi}=0.70$~rad. 
The limit on the commanded heading angle was set to be $l_m=2.5$~rad and the desired maximum tracking error was selected to be $l_e= 0.9$~rad. The model parameters were updated online twice per figure-of-eight cycle with least-squares fits to real-time measurements. The uncertainty levels of the parameters of the steering model were set to be $20\%$ of the estimated parameter values, $\delta K = 0.2\cdot K$, $\delta t_d = 0.2\cdot t_d$.


Figure~\ref{fig:trajectory1C} shows the tracking of a reference figure-of-eight path over one cycle (Flight~1). We observe that the kite position trajectory, $\xi(t):=(\theta(t),\phi(t))$, follows the reference path, $\xi^\8$. The delay compensated prediction of the position, $\xi^{t_d}(t)$, follows the reference path more closely as it is the signal used in the \ac{mpc} objective. The difference between the trajectories of the position and its prediction can be explained by model mismatch and to a greater extend by the tracking error in the inner-loop controller. 
For the same flight interval (Flight~1), the velocity perpendicular to the tethers, $v_{\theta\phi}(t)$, obtained from the derivative of the line angle measurements, and its prediction based on~\eqref{eq:velocity} are shown in Figure~\ref{fig:velocity1C}. The significant velocity variation within the cycle is captured by the velocity model in \eqref{eq:velocity} which leads to better predictions compared to the assumption of constant velocity. The tracking of the commanded heading angle during the cycle (Flight~1) is illustrated in Figure~\ref{fig:gamma1C}. The heading angle, $\gamma(t)$ follows the commanded signal, $\gamma^\re(t)$, well with a time shift. Note that the model of the steering dynamics in~\eqref{eq:steeringModel} is independent of the kite position but due to the feedback of the estimated heading angle good tracking can be achieved. 

\begin{figure}[tb]
  \centering
  \includegraphics[width=\linewidth]{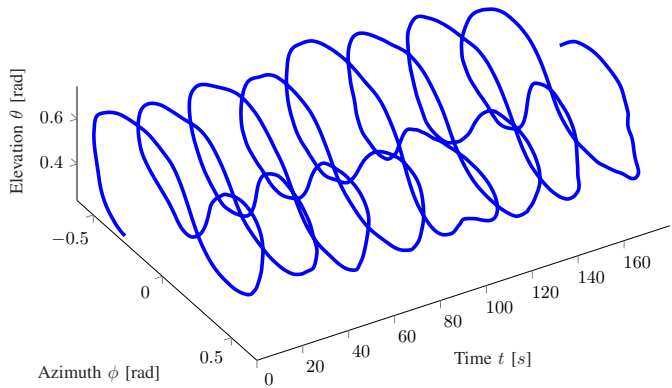}
  \caption{Long duration position trajectory, $\xi(t)$ (Flight~2).\vspace{0pt}}
  \label{fig:trajectory8C}
\end{figure}

\begin{figure}[tb] 
  \centering
  \begin{subfigure}{\linewidth}
      \includegraphics[width=\linewidth]{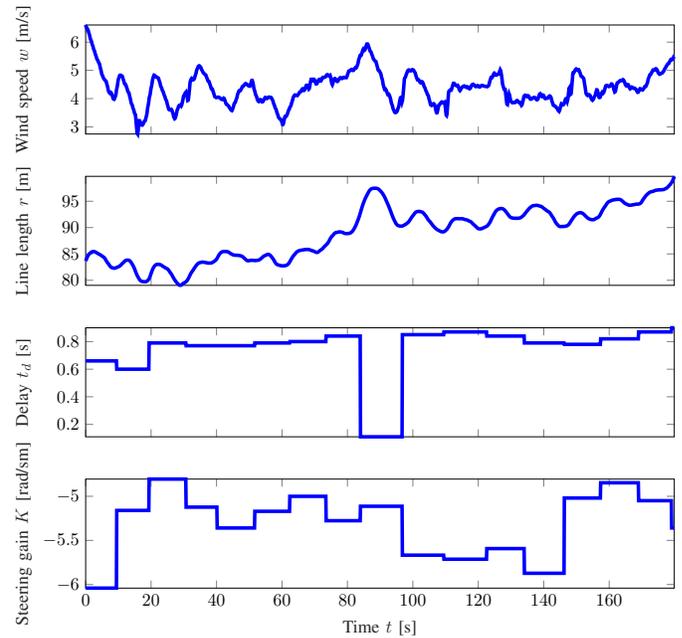}
  \caption{Measurements of the wind speed, $w(t)$, and the line length, $r(t)$; estimates of the parameters in the steering model~\eqref{eq:steeringModel}.}
  \label{fig:parameters}
  \end{subfigure}
    \begin{subfigure}{\linewidth}
      \vspace{15pt}
      \includegraphics[width=\linewidth]{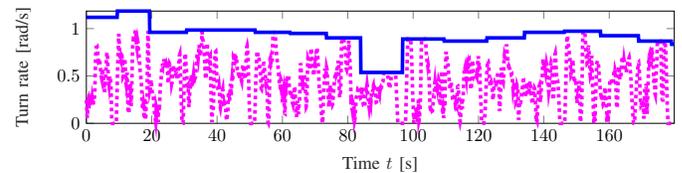}
  \caption{Bound on the rate of change of the commanded heading angle, $l_r$ ({\color{blue}solid}), and resulting magnitude  of the rate of the commanded orientation, $|\dot{\gamma}^\re(t)|$ ({\color{magenta}dotted}).}
  \label{fig:rateLimit}
  \end{subfigure}
  \caption{Evolution of parameters (Flight~2): operating conditions and parameter estimates (\ref{fig:parameters}), resulting limit on the rate of change of the commanded signal (\ref{fig:rateLimit}).\vspace{0pt}}  \label{fig:parametersRateLimit}
\end{figure}

Considering a longer flight duration (Flight~2), we can observe the adaptation to changing operating conditions. Figure~\ref{fig:trajectory8C} shows a flight trajectory over 3 minutes. 
As the wind speed and line length change, the parameter estimates vary, shown in Figure~\ref{fig:parameters}, and the bound on the rate of change of the commanded orientation, shown in Figure~\ref{fig:rateLimit}, is adapted according to the changing limitations of the tracking controller. We observe that the controller is able to track figure-of-eights cycles while the wind speed changes between 2.7~m/s and 6.6~m/s. 
Note that the periods of the reference paths are relatively high due to the large delay and the high uncertainty considered in the steering model. 

Throughout the experiment the line length is reeled out from 79~m to 100~m. For larger line lengths the estimate of the delay increases in general. For fast increases of the wind speed, however, the tether forces increase resulting in less line sag and a lower delay estimate despite the line reeling out as evident for 80-100~s in Figure~\ref{fig:parameters}. In this event the model assumptions are inaccurate leading to a poor model fit and an outlier in the parameter estimates. 

By updating the model parameters every half cycle the system can react to slow disturbances and adapt to a varying environment. The results indicate that good tracking performance can be achieved by constraining the guidance based on the limitations of the tracking controller and by re-evaluating these for changing operating conditions.
  
\section{Conclusion}\label{sec:conclusion}
We have presented an \ac{mpc} 
approach to fly kites autonomously in crosswind conditions for power generation. The control problem is split into two parts with a cascaded control architecture. We apply predictor feedback to account for the input delay affecting the system. To ensure that the commanded signal determined by the outer control loop is not too fast for the inner control loop to track, the limitations of the inner-loop controller are determined using a robustness analysis based on model parameter uncertainty. The limitations are parametrised by a bound on the 
rate of change of the commanded signal. This bound is taken into account as a constraint in the 
predictive 
guidance control. 

The main benefit of applying \ac{mpc} 
in this approach is constraint satisfaction. The optimisation framework, however, also enables the consideration of new objectives. 
The approach has been successfully tested in tow test flight experiments. The optimisation-based guidance strategy is able to steer the kite to follow figure-of-eight paths. By adapting model parameters and control constraints online using updated measurements, autonomous flights under strongly varying operating conditions were achieved. 



\section{Acknowledgements}
The authors would to thank J.\ Heilmann, D.\ Aregger, I.\ Lymperopoulos, and S.\ Diwale (A2WE project~\cite{A2WE}) for their contribution in conducting the flight experiments. Support from embotech GmbH and  in particular A.\ Hempel is gratefully acknowledged.

\bibliographystyle{IEEEtran}
\bibliography{kites,control}

\end{document}